\documentclass[aps,print]{revtex4}  

\usepackage{comment}
\usepackage{amsmath}
\usepackage{amssymb}
\usepackage{graphicx}
\usepackage{float}
\usepackage[usenames,dvipsnames]{xcolor}
\usepackage[normalem]{ulem}

\newcommand{\ham}{\mathcal{H}}

\newcommand*\hbz{ \hat{\bf z} }
\newcommand*\hbzt{{ \hat{\bf z}(t_n) }}
\newcommand*\bz{ {\bf z} }
\newcommand*\bzt{ {\bf z}(t_n) }

\newcommand*\dbz{\delta \bz}

\newcommand*\dbzdott{\dot{\delta \bz}(t_n)}
\newcommand*\dbzt{ {\delta \bz}(t_n) }
\newcommand*\Bell{\ensuremath{\boldsymbol\ell}} 
 
\newcommand*\Bellt{{\ensuremath{\boldsymbol\ell}}(t_n)}

\newcommand*\Hes{\boldsymbol{\mathcal{D}}}
\newcommand*\ecbz{{\hat{\bz}}_{ec}(t_n)}
\newcommand*\FF{{\bf F}}

\date{}
\begin{document}

\title{Dynamical Systems and Neural Networks}

\author{Akshunna S. Dogra}
\affiliation{John A. Paulson School of Engineering and Applied Sciences, Harvard University, Cambridge, Massachusetts 02138, USA}

\begin{abstract}
Neural Networks (NNs) have been identified as a potentially powerful tool in the study of complex dynamical systems. A good example is the NN differential equation (DE) solver, which provides closed form, differentiable, functional approximations for the evolution of a wide variety of dynamical systems. A major disadvantage of such NN solvers can be the amount of computational resources needed to achieve accuracy comparable to existing numerical solvers. We present new strategies for existing dynamical system NN DE solvers, making efficient use of the \textit{learnt} information, to speed up their training process, while still pursuing a completely unsupervised approach. We establish a fundamental connection between NN theory and dynamical systems theory via Koopman Operator Theory (KOT), by showing that the usual training processes for Neural Nets are fertile ground for identifying multiple Koopman operators of interest. We end by illuminating certain applications that KOT might have for NNs in general.
\end{abstract}

\flushbottom
\maketitle

\section{Introduction}
Let $\FF(\mathbf{z} = z_1, z_2,...., z_{D-1}, z_D)$ be a well behaved operator governing some $D$ dimensional dynamical flow:
\begin{equation}
    \label{eq:Dynflow}
    \dot {\bf z}  = \FF({\bf z})
\end{equation}
In \cite{mar20}, the authors presented a feed forward Neural Network (NN) that used loss functions based on Hamilton's equations to predict the evolution of various Hamiltonian systems over fixed temporal domains (that is, $\FF \equiv {\bf J}\cdot\nabla \ham$, where $\bf J$ is the symplectic matrix, $\ham$ is the system Hamiltonian and \bz is the phase space vector). The Neural Net (NN) was several orders of magnitude more accurate than a traditional symplectic Euler solver with an equivalent temporal discretization. It also provided the benefit of generating a smooth, closed form, functional approximation for the expected dynamics over the chosen temporal domain: providing precise and robust expressions for the evolution of state space parameters and thus, the entire physical system. Lastly, the NN methods could be parallelized in physical time, something iterative methods are incapable of providing. 

The chief disadvantage of the NN, when compared to existing methods, was the amount of computational resources required to generate a single solution. These computational costs severely limited the utility of the particular NN solver, especially in a low dimensional physical system or when domain resolution was not expected to be of substantial importance. Further, identifying the local error in NN prediction was analytically intractable, unlike the case for traditional numerical solvers, where one may occasionally compute exact local error terms to refine the prediction further. 

A plurality of physical systems are described by state variables that flow smoothly in time, save perhaps for a finite number of discontinuities. In this work, we extend the applicability of feed forward NNs similar to the one presented in \cite{mar20} beyond Hamiltonian systems, working with the weaker assumption that $\FF$ is simply a smooth operator over the domain of interest. Most of the strategies we describe are also valid for operators $\FF$ that are smooth globally, except at a finite number of simply connected regions, but we do not present formal proofs for the same. We derive results about the error in the NN predictions and prescribe error correction methods to magnify the speed and accuracy of the NN. We demonstrate how the training process is intimately connected with Koopman Operator Theory, allowing us to describe a series of new strategies for NN training, including some that provide savings in computational costs of training NNs. We end on a discussion of the challenges the undertaking described above would face and how those challenges could be overcome to make \textit{Koopman training} a viable technique for training NNs.

\section{Neural Network Approaches to Smooth dynamical systems}
The central aim of this section is to generalize existing results and showcase a strategy that can aid the training process of existing NNs like the ones described in \cite{lag98}\cite{mar20}, without modifying architectures and minimal additional computational resources.

In \cite{mar20}, the authors presented a rapidly convergent NN that could find accurate functional approximations $\hbz(t)$ for the evolution of phase space parameters $\bz(t)$ of various Hamiltonian systems - chaotic and nonlinear systems included - by simply demanding information about the initial phase space co-ordinate $\bz(0)$ and the temporal domain $[0,T]$ of interest. Through a series of examples and results, the authors demonstrated how NN training can both be optimized by tailoring the architecture to the problem at hand (\textit{physical insight} optimizing NN operations by advising the choice of activation functions) and how physical parameters of inherent significance can be studied better by involving the recent advances that machine learning methods, specifically deep NNs, have made in the past few decades (NNs bettering \textit{physical insight} by providing accurate approximations to the dynamics at hand). 

The NN itself was structurally simple: an input layer demanding a set of $N$ randomly generated points in the temporal domain of interest $[0,T]$ each iteration, two hidden layers with $sin()$ activation hubs and an output layer with $D$ outputs ${\bf N} \equiv \{N_1, N_2, ..., N_D\} $ - one for each state parameter described in Eqn. \ref{eq:Dynflow}. Thus, the NN was a $D$ - dimensional output map for $t$: sourcing $t$ from the temporal domain of interest meant the NN was being trained to be an effective functional approximation $\hbz(t)$ for the expected evolution $\bz(t)$ of the dynamical system over the chosen temporal domain. To enforce the initial condition during the training process, the NN output hubs ${\bf N}(t)$ and the final NN prediction $\hbzt$ were related as: $\hbzt = \bz(0) + (1-e^{-t}){\bf N}(t)$.

The authors used the symplectic form of Hamilton's equations as the basis for the mean squared temporal loss. 
\begin{equation}
    \label{mar20_loss}
    L = \overline{\Bell(t_n) \cdot \Bell(t_n)} \text{ :}\text{ }\text{ }\text{ }\text{ }\text{ }\text{ }\text{}\text{ }\text{ }\text{ }\text{ }\text{ } \Bell (t_n) = {\dot{\hbz}}(t_n) - {\bf J} \cdot (\nabla \ham)|_{\hbzt}
\end{equation}
where $\bf J$ is the symplectic matrix and $(\nabla \ham)|_{\hbzt}$ is the gradient of the system Hamiltonian, evaluated at the prediction $\hbzt$. Since $\hbz$ was a function of $t$ by the construction of the network, $\dot\hbz$ could be evaluated at any $t_n$ from within the NN, making the training completely unsupervised. 

Eqn. \ref{mar20_loss} also implies the capacity of the NN to dispense with causality, since the evaluation of $\Bellt{}$ is not dependent on the evaluation of $\Bell (t_{n-1})$. This implies that at least part of NN training could be parallelized.
The authors also presented the following repackaging of the loss function:
\begin{equation}
    \label{eq:Err_anal_1}
    L = \overline{\Bell(t_n) \cdot \Bell(t_n)}\text{ :}\text{ }\text{ }\text{ }\text{ }\text{ }\text{ }\text{}\text{ }\text{ }\text{ }\text{ }\text{ } \Bell (t_n) \approx   {\bf J} \cdot (\Hes(\ham)|_{\hbz (t_n)}\cdot \delta{\bz}(t_n))   -  \dot{\delta \bz}(t_n)
\end{equation}
where $\Hes(\ham)|_{\hbz (t_n)}$ is the Hessian matrix of the system Hamiltonian evaluated at $\hbzt$ and $\delta\bz (t_n)= \bz (t_n) - \hbz (t_n)$ is the difference between the true evolution and the NN prediction.

Let us describe (and generalize to smooth dynamical systems) strategies for speeding the training of NN architectures similar to the ones found in \cite{lag98}\cite{mar20}. Let the NN make its predictions for some discrete, finite set of $N+1$ time points $\{t_n \}$, with $t_0=0$ and $t_{N}=T$ being the endpoints of our temporal domain of interest. All intermediate $t_n$ are sampled randomly from the set $(0,T)$ before each forward pass. We define the following time averaged function $L$ as our \textit{cost/loss} function:
\begin{equation}
    \label{eq:Err_anal_initial}
    L = \overline{\Bell(t_n) \cdot \Bell(t_n)}\text{ :}\text{ }\text{ }\text{ }\text{ }\text{ }\text{ }\text{}\text{ }\text{ }\text{ }\text{ }\text{ } \Bell (t_n) = {\dot{\hbz}}(t_n) - \FF(\hbzt)
\end{equation}
where $\hbzt$ is the predicted output and $\FF$ is the dynamical operator describing the system. Let $\bzt$ be the true value at $t_n$ and $\dbzt = \bzt - \hbzt$. Let us assume the network is trained sufficiently such that the Taylor expansion is convergent throughout the temporal domain. We get:
\begin{equation}
\FF(\bz) = \FF(\hbz) + (\FF _{\bz}|_{\hbz}\cdot\dbz) + (\dbz ^T\cdot\FF _{\bz\bz}|_{\hbz}\cdot\dbz) + ... 
\end{equation}
\begin{equation}
\label{eq:Taylor_initial}
\implies \FF(\hbz) = \FF(\bz) - [(\FF _{\bz}|_{\hbz}\cdot\dbz) + (\dbz ^T\cdot\FF _{\bz\bz}|_{\hbz}\cdot\dbz) + ...]
\end{equation}
Here, $\FF_{\bz}|_{\hbz}\equiv \nabla \FF$, $\FF_{\bz\bz}|_{\hbz}\equiv \nabla(\nabla \FF)$ and so on, evaluated at $\hbzt$. We note that many common dynamical operators are built from elementary functions with infinite or qualitatively large radii of convergence. 

We know that $\forall t_n, {\dot{\bz}}(t_n) - \FF(\bzt) = 0$. We use Eqn. \ref{eq:Err_anal_initial} and \ref{eq:Taylor_initial} to generalize Eqn. \ref{eq:Err_anal_1} (Eqn. 15 in \cite{mar20}):
\begin{equation}
    \label{eq:Err_anal_reworked}
    L = \overline{\Bell(t_n) \cdot \Bell(t_n)}\text{ :}\text{ }\text{ }\text{ }\text{ }\text{ }\text{ }\text{}\text{ }\text{ }\text{ }\text{ }\text{ } \Bell (t_n) = [(\FF _{\bz}|_{\hbzt}\cdot\dbzt) + (\dbz^T(t_n)\cdot\FF _{\bz\bz}|_{\hbzt}\cdot\dbzt) + ...] - \dbzdott
\end{equation}
Let us say that for some given loss profile, $|\Bell (t)| \leq \ell _{max}$ and $\sigma _{min}$ is the minimum singular value of $\FF _{\bz}$ over $[0,T]$. Keeping only the leading order $\FF _{\bz}|_{\hbzt}\cdot\dbzt$ term in Eqn. \ref{eq:Err_anal_reworked} and using the exact methodology as the authors in \cite{mar20}, we obtain the generalization of their result (Eqn. $24$ in \cite{mar20}) on the upper bounds on magnitudes of individual error components of $\dbz$ as:
\begin{equation}
\label{generalized_error_bounds}
\|\delta z_{i}\| \leq \frac{\ell_{\text{max}}}{\sigma_\text{min}} 
\end{equation}
Eqn. \ref{generalized_error_bounds} gives us a natural way of using the loss function to bound the error in the predicted solution by the NN. However, given smoothness for our operator $\FF$ and the sufficient training assumption used to derive Eqn. \ref{eq:Err_anal_reworked}, we can derive more than just bounds on the error we should expect from our NN.

Eqn. \ref{eq:Err_anal_initial} tells us $\Bell(t)$ is a smooth function (since $\hbz$ is a smooth function of $t$ by the construction of our NN and the operator $\FF$ is assumed to be a smooth operator). Further, the NN can calculate $\Bellt, \FF _{\bz}|_{\hbzt}, \FF _{\bz\bz}|_{\hbzt}, ...$ for any $t_n$.
Therefore, Eqn. \ref{eq:Err_anal_reworked} allows us the following discrete, recursive equation to estimate $\dbzt$ by picking a small enough $\Delta t_n$, such that it is reasonably resolved (such a finite $\Delta t_n$ exists due to smoothness of $\Bell$ and \FF):
\begin{equation}
    \label{Err_anal_fixing_term}
    \dbz(t_{n+1}) = \dbzt + \Delta t_n[[(\FF _{\bz}|_{\hbzt}\cdot\dbzt) + (\dbz^T(t_n)\cdot\FF _{\bz\bz}|_{\hbzt}\cdot\dbzt) + ...] - \Bellt] : \text{ }\text{ }\text{ } \Delta t_n = t_{n+1} - t_n,\text{ } \dbz (0) = 0
\end{equation}
We can use Eqn. \ref{Err_anal_fixing_term} to generate error data about the NN prediction to as good a resolution and accuracy as needed, by choosing an adequately small $\overline{\Delta t_n}$. This capability is quite useful.

One of the major computational costs of training the NNs in \cite{lag98}\cite{mar20} is the calculation of $\dot{\hbz}$. This is because $L$ is defined as the mean squared residual of the governing dynamical equation (Eqn. \ref{eq:Err_anal_initial}): hence, the differential equation has to be solved before backpropagation can be applied. An efficient way of cutting down this cost would be to:
\begin{enumerate}
    \item train the network until $\hbzt$ is reasonably within the radius of convergence of $\bzt$ for all $t_n$, as assumed in Eqn. \ref{eq:Err_anal_reworked} (Eqn. \ref{generalized_error_bounds} can identify $L$ at which that may reasonably be assumed)
    \item generate $\dbzt$ for the NN at a precise enough resolution at the end of the $k^{th}$ training iteration
    \item produce an error corrected prediction dataset $\ecbz=\hbzt + \dbzt$
    \item redefine $\Bellt: \Bellt = \hbzt - \ecbz$
    \item keep $t_0 = 0, t_N=T$ and select $N-1$ other time points to assemble batches randomly from the set of $t_n$ for which $\bz _{ec}(t_n)$ is available (generating the error data at a higher resolution also provides the benefit that the likelihood of repeating batches is very low, mimicking SGD)
    \item train the NN, using the new $\Bellt$ definition, until $L$ is small enough that $\ecbz$ needs to be more accurate. 
    \item redefine $\Bellt: \Bellt = {\dot{\hbz}}(t_n) - \FF(\hbzt)$. Train the NN using original setup for a few iterations. Repeat $1-6$.
\end{enumerate}
This algorithm cuts down the computational costs of calculating ${\dot{\hbz}}(t_n)$ (major NN computational cost) and $\FF(\hbzt)$ (a minor cost) for each future iteration. For NNs described in \cite{lag98}\cite{mar20}, this is the dominant computational cost per iteration. Let us assume the precision needed to generate sufficiently resolved $\dbzt$ is $k$ times higher than the temporal discretization being used per iteration in the original NNs. Then, one would consume at worst $k$ forward pass computational cost to setup the process. However, simple combinatorics dictates that the NN could train using the redefined $\Bellt{}$ for practically any number of iterations if either $k$ and/or the batch size was large enough. To put this into perspective, for the NN presented in \cite{mar20}, $k=2$ would practically ensure that an exactly repeated batch never occurs with random selection. In the same NN, the utility of $k=10$, in terms of saved costs, would last orders of magnitude more iterations longer than the $10$ forward passes needed to setup the process, before over-fitting concerns start building up appreciably.

\section{Koopman Operators for Training a Neural Net}
Recent work by Redman \cite{WTR} has shown that the Renormalization Group (RG), a powerful tool in theoretical physics, is intimately connected with Koopman Operator Theory (KOT), a sub-field of dynamical systems theory. By presenting the block spin renormalization process as a dynamical flow in the space of coupling constants and iteration step as the temporal parameter, Redman proved that the RG is a Koopman Operator (KO) by definition. This realization was used in conjunction with algorithms inspired by KOT to compute useful information about critical exponents of different physical systems, without the translational invariance assumption that drastically limits the applicability of RG theory to complex physical systems.

We provide similar results connecting KOT with NNs, by repackaging NN training as a discrete dynamical flow, with the number of iterations serving as the discrete temporal parameter. Some dynamical quantities of interest in this flow are the loss function $L$ (or its independent components) and the individual weights of the NN. 

From a purely computational/technical standpoint, NN training is done to find a set of weights that minimize the loss $L$ to the best possible extent (while still serving the interests of the model being built. However, that detail can be considered extraneous to the dynamical picture we are trying to build and profit from in this work). From a dynamical systems standpoint, the object is to find the stable fixed point(s) for $L$ or loss function components $L_1, L_2, ..., L_k$ (state function(s) of interest), whose evolution is governed by the evolution of different weights (state space parameters) as the number of training iterations increases (temporal parameter for the system). This reasoning can also be flipped to envision a complementary dynamical setup for the same NN architecture and training process - the need to find stable fixed point(s) for each weight (state function of interest), whose evolution is governed by the loss function $L$ or its components (state space variables). In this picture, we track the evolution of individual weights as a result of the evolution of the NN loss $L$ (or loss components $L_1,L_2,...$). 
We will show the existence of one KO of interest for each of the descriptions above. These results, paired with existing KOT algorithms to efficiently construct and use Koopman Operators, predict new tools to optimize and benchmark the training of a NN.

Let us give a short description of KOT. It was created in 1931 by Bernard Koopman and later expanded upon by Koopman and Neumann in 1932 \cite{KOT_32}. KOT provides a spectral approach to dynamical systems, specially in the context of nonlinear systems, by investigating the spectrum of the Koopman (or composition) operator $U$ of state space functions $f({\bf w}\equiv \{w_1,w_2, ..., w_M\})$, $\bf w$ being the state space parameters. Here, the KO $U$ is the infinite dimensional linear operator satisfying:
\begin{equation}
    \label{KOT_defn}
    U^tf({\bf w_0}) = f(g^t(\bf w_0))
\end{equation}
where $g^t$ prescribes the evolution of the state parameters. KOT has many interesting facets, but for the purposes of this work, we shall focus on two: its capacity to supply the dynamics of state functions of interest and its capacity to identify positive invariant sets and limit points of dynamical flows - especially when such a search is driven by data \cite{IM_05}. In particular, we shall discuss the possibility of data driven KOs replacing the \textit{standard training} of a NN and/or identifying the limits to which a NN can be trained and the corresponding weight values that should be associated with those limits (or the limit point(s) of the weight flow and the corresponding $L$ (or $L_1, L_2,...$) associated with those limit point(s)).

Let $L({\bf w})$ be a scalar loss function for some NN, where ${\bf w}\equiv \{w_1,w_2, ..., w_M\}$ represents the $M$ weights of a NN in some order. A successful training process updates the weights at each backward pass to minimize the loss function, seeking a possible configuration of weights such that $L = 0$ (or tends to some other local minimum value, depending upon how the loss function is defined and the efficacy of the training procedure). For the purposes of this discussion, we are dispensing with caveat that a lower $L$ might not necessarily mean a better NN (over-fitting being the most obvious pitfall). We are simply interested in whether there is merit to the notion that the standard training process could be replaced by its \textit{Koopman} analogue or not. Put another way, we are interested in mimicking the dynamical flow caused by the standard training methods using KOT tools.

Thus, the training rule $U$ for a NN is a discrete dynamical map on $L$, with the training iteration $t$ being the temporal parameter. Let $T$ be the discrete mapping governing the dynamics for the weights. Then,
\begin{equation}
    \label{Loss_KOT}
    L(t) = U^tL({\bf w_0}) = L(T^t({\bf w_0}))
\end{equation}
where ${\bf w_0}$ describes some \textit{original} set of weights for our NN and $t$ is the number of iterations from our temporal origin. The second equality shows that the training $U$ is, by definition, a Koopman operator on the state function $L$.

The linearity of the KO is evident from Eqn. \ref{Loss_KOT}. For example, many loss functions are themselves the sum of constituent, independent loss function like components $L_1, L_2, ...$, such as those in \cite{lag98}\cite{mar20}.
\begin{equation}
    \label{Linear_Loss_KOT}
    U^t[L_1({\bf w_0})+L_2({\bf w_0})+...] = [L_1+L_2+...](T^t({\bf w_0})) = L_1(T^t({\bf w_0})) + L_2 (T^t({\bf w_0})) + ... = U^t L_1({\bf w_0}) + U^t L_2({\bf w_0}) + ...
\end{equation}
For the kind of NNs discussed until now, another KO with possible applications is identifiable. Let $w$ be some arbitrary weight of the NN. The value of $w$ at the end of each iteration is also governed by some training rule $U$. Let $U$ now represent that discrete map for $w$ and $T$ be the discrete map supplying the dynamics for the loss components  $\{L_1, L_2, ...\} = \bf L$. We have:
\begin{equation}
    \label{weight_KOT}
    w(t) = U^tw({\bf L_0}) = w(T^t({\bf L_0}))
\end{equation}
where   $\{L_1, L_2, ...\} = \bf L_0$ describes the \textit{original} value of loss components.
Eqn. \ref{Loss_KOT} and \ref{weight_KOT} clearly demonstrate that the training of NNs imposes a dynamical flow on various inherent quantities of interest and the mapping associated with each training is a Koopman operator in this dynamical picture.

Having established the basic foundation, what benefits can Eqn. \ref{Loss_KOT} and \ref{weight_KOT} provide? Eqn. \ref{Loss_KOT}, \ref{weight_KOT} imply that if we could find the associated KOs with their respective flows, we would obtain new tools for both training NNs and identifying their capabilities.  Assuming that an adequate approximation $\hat{U}$ to the operator $U$ can be found, Eqn. \ref{Loss_KOT} and \ref{weight_KOT} can update the relevant state functions (whether weights or the loss) without having to go through the standard forward-backward propagation loop. If $L$ (or $L_1, L_2, ...$) is (are) our state function, we can use the associated $\hat{U}$ to figure out the limits of a particular architecture and training, without having to do the complete training itself. If the state functions are the weights, we can use the associated $\hat{U}$ as a substitute for the forward-backward propagation based training, potentially saving immense computational resources. We call this way of updating weights \textit{Koopman training} and it should be expected to be especially powerful when the weights (or loss function) enter the attraction basin of any of their fixed points.

Let us say some fully connected, feed forward NN has $m$ layers, each with $m$ neurons. We intend to create approximations $\hat{U}$ to KOs that governs the flow of weights connecting each neuron of a preceding layer to some arbitrary neuron in the next one (this would be akin to building a separate KO for each column in the connection matrix listing weights connecting neurons of two successive layers. There are computationally more efficient choices available, but we are going with a simple case). We assume that we will be using the Frobenius–Perron operator inspired methods detailed in \cite{IM_19} to construct the KO. Finally, let us assume that we need data from $k$ iterations to build approximate KOs. The total construction cost of \textit{all} the relevant $\hat{U}$ is bounded above by order $km^4$. The total per iteration cost is bounded above at order $m^4$ (since Koopman training simply involves multiplying $m^2$ matrices ($m \times m$ type) with $m^2$ vectors (dimension $m$)). Hence, even the brute force approaches to \textit{Koopman training} have the promise of computational savings, since creation and usage of KOs do not require any calculations beyond those incurred during direct arithmetical matrix multiplications. The lack of intermediate steps every iteration means KOs don't have large constant prefactors in their usage complexity. Contrast this with the NN presented in \cite{lag98}\cite{mar20}, where the calculation of derivative terms for the evaluation of the loss function alone is a substantial additional computational cost. Lastly, the major sink for computational costs when training a NN is at the latter half of training, when the training provides diminished returns per iteration, due to the NN having found some sort of a stable regime for the weights or $L$ - exactly the kind of regime where linear approaches to non-linear phenomena come in especially handy and is extensively studied. 

Let us say we somehow know spectral information about the operator $U$. We show the utility of such knowledge by assuming we are investigating one of the simpler kinds of dynamical systems. We define a Koopman eigenfunction-eigenvalue pair as objects that satisfy the following:
\begin{equation}
    \label{KOT_eigen}
    U\phi _i = e^{\lambda _i}\phi _i \implies U^t\phi _i = e^{\lambda _i t}\phi _i \text{ :}\text{ }\text{ }\text{ }\text{ }\text{ }\text{ }\text{ }\text{ }\text{ }\text{ }\text{ }\text{ }\text{ }\text{ }\text{ }\text{ }\text{ } \lambda _i \in \mathbb{C}
\end{equation}
For a wide class of well behaved linear and nonlinear dynamical systems \cite{IM_19}, state functions lie in the span of such eigenfunctions:
\begin{equation}
    \label{KOT_eigen_decomp}
    w = \sum_{i=1}^{\infty}w_i \phi _i \text{ }\text{ }\text{ }\text{ or}\text{ }\text{ }\text{ }\text{ }L = \sum_{i=1}^{\infty}L_i\phi _i
\end{equation}
For such systems, the action of $U$ on the state function $w$ is equivalent to the following expression:
\begin{equation}
    \label{KOT_Decomp}
    U(w) = \sum_{i=1}^{\infty}w_i e^{\lambda _i} \phi _i \implies U^t(w) = \sum_{i=1}^{\infty}w_i e^{\lambda _i t} \phi _i \text{ }\text{ }\text{ }\text{ or}\text{ }\text{ }\text{ }\text{ }U(L) = \sum_{i=1}^{\infty}L_i e^{\lambda _i} \phi _i \implies U^t(L) = \sum_{i=1}^{\infty}L_i e^{\lambda _i t} \phi _i
\end{equation}
where $\phi _i$ are the Koopman eigenfunctions and $L_i$ are the Koopman modes. 

Eqn. \ref{KOT_Decomp} has even more powerful implications - if the spectral information about the KO associated with a particular weight (or loss) flow is estimable, then the computational complexity of \textit{Koopman training} is reduced practically to that of finding the leading Koopman modes $L_i$ and eigenfunctions $\phi _i$. Koopman modes also make it possible to recover the evolution of state space parameters themselves, by studying the action of the KO on the identity state function. Thus, a single \textit{Koopman training} approach to NNs can also reveal limit points for both the loss components and the weights. Unfortunately, identifying KO spectral features analytically for even simple nonlinear systems is often a herculean task. Fortunately, data provides another way. 

A major factor in the resurgence of KOT has been the work of Mezić \cite{IM_19}\cite{IM_05}\cite{IM_19_2} and others\cite{MO_15}\cite{QL_17}, who have advanced the field by introducing data driven methods to efficiently use KOT over a wide variety of deterministic and stochastic dynamical systems. In particular, the major achievements of data driven KOT has been its capacity to identify interesting objects like invariant sets and limit points for a wide variety of nonlinear systems - applications that are directly relevant and applicable to NNs. 

The construction details for approximations $\hat{U}$ to KOs associated with a large variety of arbitrary flows, both deterministic and stochastic, can be found in \cite{IM_19}. Data driven methods for Koopman mode decomposition were derived and demonstrated in \cite{MO_15}. NNs themselves have been used to optimize the identification of the leading Koopman modes, eigenfunctions and eigenvalues \cite{QL_17} - a potential setup for an approach where NNs help train NNs (or conversely, KOs help identify KOs). Thus, there exist clear, proven methods for constructing \textit{Koopman training} analogs for a large variety of standard NN training methods (gradient and stochastic gradient descent (GD/SGD) included). We hope that the realization that NNs themselves belong to the class of dynamical flows that are accessible to data driven KOT will encourage future developments in both machine learning and KOT.

\begin{acknowledgments}
We acknowledge the financial support provided by the School of Engineering and Applied Sciences at Harvard University, with a special thanks to its members: Dr. Marios Mattheakis, Dr. David Sondak, Dr. Pavlos Protopapas and Prof. Efthimios Kaxiras. We thank Mr. William T. Redman for discussions on Koopman Operators.
\end{acknowledgments}

\end{document}